\documentclass[12pt,reqno]{amsart}

\usepackage{amssymb,latexsym}

\usepackage{enumerate}
\allowdisplaybreaks
\usepackage[french,english]{babel}
\usepackage{amsmath}
\usepackage{graphicx}
\usepackage{amssymb}
\usepackage{bbm}
\usepackage{amsthm,mathtools}
\usepackage{ulem}
\usepackage{geometry}
\usepackage{tikz-cd}
\usepackage{mathrsfs}
\usepackage[colorinlistoftodos]{todonotes}
\usepackage{enumitem}
\usepackage{verbatim}
\usepackage[foot]{amsaddr}
\usepackage{dsfont}
\usepackage{cite}
\usepackage[T1]{fontenc}

\makeatletter

\@namedef{subjclassname@2010}{
	
	\textup{2020} Mathematics Subject Classification}

\makeatother
\newtheorem{thm}{Theorem}[section]
\newtheorem*{thm*}{Theorem}

\newtheorem{lem}[thm]{Lemma}
\newtheorem{pro}[thm]{Proposition}
\theoremstyle{definition}

\numberwithin{equation}{section}

\newcommand{\kf}{\kappa_\mathbb{F}}
\newcommand{\zf}{\zeta_\mathbb{F}}
\newcommand{\af}{a_\mathbb{F}}
\newcommand{\bg}{\big}

\newcommand{\Bg}{\Big}

\newcommand{\lo}{\log_2}
\newcommand{\lt}{\log_3}
\newcommand{\lf}{\log_4}
\newcommand{\inv}{^{-1}}

\newcommand{\mbc}{\mathbb{C}}

\newcommand{\mbe}{\mathbb{E}}

\newcommand{\mbr}{\mathbb{R}}

\newcommand{\mbf}{\mathbb{F}}

\newcommand{\mbn}{\mathbb{N}}

\newcommand{\mbq}{\mathbb{Q}}

\newcommand{\mch}{\mathcal{H}}

\newcommand{\mcb}{\mathcal{B}}

\newcommand{\mce}{\mathcal{E}}

\newcommand{\mcp}{\mathcal{P}}
\newcommand{\mck}{\mathcal{K}}

\newcommand{\mco}{\mathcal{O}}
\newcommand{\mcq}{\mathcal{Q}}
\newcommand{\mcm}{\mathcal{M}}
\newcommand{\mcr}{\mathcal{R}}

\newcommand{\mmd}{\mathrm{d}}
\newcommand{\mme}{\mathrm{e}}
\newcommand{\mmi}{\mathrm{i}}
\newcommand{\mit}{\mathrm{i}t}
\newcommand{\mdt}{\mathrm{d}t}
\newcommand{\re}{\operatorname{Re}}

\usepackage{hyperref}
\hypersetup{colorlinks=true,linkcolor=blue,anchorcolor=blue,citecolor=blue}
\usepackage{color}
\newcommand{\newabstract}[1]{%
	\par\bigskip
	\csname otherlanguage*\endcsname{#1}%
	\csname captions#1\endcsname
	\item[\hskip\labelsep\scshape\abstractname.]
}

\begin{document}

	\baselineskip=17pt

	\title[Z. Dong, Q. Yang, S. Zhao and G.-L. Zhou]{Extreme values of derivatives of Dedekind zeta functions on the \(1\)-line}
	
    \author{Zikang Dong\textsuperscript{1}}
	\author{Qiyu Yang\textsuperscript{2}}
    \author{Shengbo Zhao\textsuperscript{3}}
    \author{Guang-Liang Zhou\textsuperscript{4}}

    \address{1. School of Mathematical Sciences, Soochow University, Suzhou 215006, P. R. China}
	\address{2. School of Mathematics and Statistics, Henan Normal University, Xinxiang 453007, CHINA}
	\address{3. School of Mathematical Sciences, Key Laboratory of Intelligent Computing and Applications (Tongji University), Ministry of Education, Tongji University, Shanghai 200092, China}
    \address{4. Department of Applied Mathematics, Nanjing Forestry University, Nanjing 210037, People's Republic of China}

    \email{zikangdong@gmail.com}
	\email{qyyang.must@gmail.com}
	\email{shengbozhao@hotmail.com}
    \email{guangliangzhou@126.com}

	\begin{abstract} 
	   In this paper, we establish lower bounds for extreme values of derivatives of Dedekind zeta functions on the \(1\)-line over arbitrary fixed number fields. Our result extends Yang's corresponding estimate for the Riemann zeta function in 2024, with the leading constant involving residues of Dedekind zeta functions and a moment of the Dickman function. A key ingredient is a uniform asymptotic formula for weighted friable sums involving the coefficients of Dedekind zeta functions.
	\end{abstract}
	
    \keywords{Extreme values, Dedekind zeta functions, resonance method, smooth numbers}
	
	\subjclass[2020]{Primary 11M41, 11M06, 11N37, 11R42.}
	
	\maketitle

\section{Introduction}
\label{sec-int}

The study of extreme values of the Riemann zeta function \(\zeta(s)\) has a long history. The \(1\)-line is the boundary of absolute convergence of its Dirichlet series and Euler product, and larger values on this line reflect the combined contributions of primes. Understanding these values also provides a model for studying the distribution of values of more general zeta and $L$-functions.

Granville and Soundararajan \cite{Granville2006RamanujanMS} proved that there
exist arbitrarily large $t$ such that
\[|\zeta(1+\mit)| \ge \mme^\gamma(\lo t+\lt t-\lf t-C),\]
where $C>0$ is an absolute constant and $\gamma$ denotes the Euler-Mascheroni constant. Here and throughout this paper, we write $\log_j$ for the $j$-th iterated logarithm. They also conjecture about the maximum of \(|\zeta(1+\mit)|\). Subsequently, Aistleitner, Mahatab and Munsch \cite{Aistleitner2019IMRN} improved the above lower bound by removing the $\lf t$ term, up to an absolute constant.

We now turn to Dedekind zeta functions associated with number fields. Let \(\mbf/\mbq\) be a fixed number field of degree \(d\), and let \(\zf(s)\) denote its Dedekind zeta function. We also write 
\[\kf=\operatorname*{Res}_{s=1}\zf(s).\]
In 2021, Dixit and Mahatab \cite[Corollary 1.4]{Dixit2021BAustMS} established a lower bound for extreme values of Dedekind zeta functions on the \(1\)-line. More precisely, there exists a constant \(C_\mbf>0\) such that
\[\max_{t \in [\sqrt{T},T]}|\zf(1+\mit)|\ge\mme^\gamma \kf(\lo T+\lt T-C_\mbf)\]
for all sufficiently large $T$. Away from the $1$-line, Balakrishnan \cite{Balakrishnan1986AA} obtained extreme values of Dedekind zeta functions of quadratic fields in the critical strip. More recently, Bondarenko, Darbar, Hagen, Heap and Seip \cite{Bondarenko2023BLMS} established lower bounds for cyclotomic fields on the critical line.

Turning to derivatives, Yang \cite{Ddyang2022Mathematika} established lower bounds for extreme values of derivatives of the Riemann zeta function, which were subsequently improved on the $1$-line by Dong and Wei \cite{Zkdong2023BAustMS}. Let $\rho$ denote the Dickman function, and put
\[
Y_j:=\int_0^\infty u^j\rho(u)du.
\]
More recently, Yang \cite[Theorem 1]{Ddyang2024BLMS} proved that, as $T\to\infty$, uniformly for all positive integers $\ell\le \log_3 T/\log_4 T$, we have
\[\max_{t \in [T,2T]}\bg|\zeta^{(\ell)}(1+\mit)\bg| \ge \bg(Y_\ell+o(1)\bg)(\lo T)^{\ell+1}.\]
For related results concerning derivatives of the Dedekind zeta function of a cyclotomic field on and near the critical line, see \cite{Zhli2026JAustMS}. Further results on extreme values of derivatives of the Riemann zeta function and Dirichlet \(L\)-functions can be found in \cite{Xhcui2026BAustMS,Zkdong2025RamanujanJ,Kalmynin2018AA,Qyyang2026JMAA}.

In the present paper, we extend Yang's result on the $1$-line to Dedekind zeta functions of arbitrary fixed number fields. Our main result is as follows.

\begin{thm}
    \label{thm1}
    Let \(\mbf/\mbq\) be a fixed number field. Then, as \(T\to\infty\), uniformly for all positive integers \(\ell \le (\lt T)/(\lf T)\), we have
    \[\max_{t \in [T,2T]}\bg|\zf^{(\ell)}(1+\mit)\bg| \ge \bg(\kf Y_\ell+o_\mbf(1)\bg)(\lo T)^{\ell+1}.\]
\end{thm}

To explain the main difficulty in extending Yang's result, for each positive integer \(n\), define 
\[
a_\mbf(n):=\#\{ \mathfrak a\subseteq\mco_\mbf: N\mathfrak a=n\},
\]
where $\mco_\mbf$ denotes the ring of integers of $\mbf$ and $N\mathfrak a$ is the absolute norm of the integral ideal $\mathfrak a$. Then, for \(\re(s)>1\), we have
\[\zf(s)=\sum_{n=1}^\infty \frac{\af(n)}{n^s}\quad\text{and}\quad(-1)^\ell\zf^{(\ell)}(s) =\sum_{n=1}^{\infty}\frac{\af(n)(\log n)^\ell}{n^s}.\]
Although \(\af(n)\) is multiplicative, the coefficients \(\af(n)(\log n)^\ell\) are not so far as \(\ell \ge 1\). Consequently, the Euler product argument underlying the long resonance method of Aistleitner, Mahatab and Munsch \cite{Aistleitner2019IMRN} cannot be directly applied to derivatives. Instead, we employ Soundararajan's method \cite{Soundararajan2008MathAnn} with a resonator supported on the divisors of a suitably chosen integer. The resulting lower bound is reduced to the estimation of weighted friable sums involving \(\af(n)\), for which we establish a uniform asymptotic formula.

Although Theorem \ref{thm1} is stated for positive integers \(\ell\), the same argument also applies to \(\ell=0\). Since
\[Y_0 = \int_0^\infty \rho(u)\mmd u =\mme^\gamma,\]
we obtain
\[\max_{t \in [T,2T]}|\zf(1+\mit)| \ge (\mme^\gamma\kf+o_\mbf(1))\lo T.\]
This agrees with the leading term of the result of Dixit and Mahatab \cite{Dixit2021BAustMS}, whose lower bound also contains a secondary term of order \(\lt T\), albeit over a different range of \(t\).

The range of \(\ell\) in Theorem \ref{thm1} can be slightly extended. In fact, for any fixed \(0<\theta<1\), the same conclusion holds uniformly for 
\[1 \le \ell \le \frac{\lt T}{\lf T-\theta \log_5 T}.\]
Indeed, the arguments of Lemmas \ref{lem-basic} and \ref{lem-tail}, together
with Stirling's formula, remain valid in this range. We retain the simpler bound in Theorem \ref{thm1} for convenience.

Theorem \ref{thm1} also yields several consequences for special classes of number fields. When \(\mbf=\mbq\), we have \(\kf=1\) and hence recover \cite[Theorem 1]{Ddyang2024BLMS}. For a fixed quadratic field \(\mbf=\mbq(\sqrt{D})\) with fundamental discriminant $D$, let \(\chi_D\) denote the corresponding quadratic Dirichlet character. The factorization
\[\zeta_{\mbq(\sqrt{D})}(s)=\zeta(s)L(s,\chi_D)\]
gives \(\kf=L(1,\chi_D)\). Consequently,
\[\max_{t \in [T,2T]}\bg|\zeta_{\mbq(\sqrt{D})}^{(\ell)}(1+\mit)\bg| \ge \bg(L(1,\chi_D)Y_\ell+o_D(1)\bg)(\lo T)^{\ell+1}.\]

The arithmetic meaning of $L(1,\chi_D)$ is given by the analytic class number formula. Let $h_D$ denote the class number of $\mbq(\sqrt D)$. Then, 
\[
L(1,\chi_D)=
\begin{cases}
\displaystyle \frac{2\pi h_D}{w_D\sqrt{|D|}}, & D<0,\\[6pt]
\displaystyle \frac{2h_D R_D}{\sqrt D}, & D>0,
\end{cases}
\]
where $w_D$ is the number of roots of unity and $R_D$ is
the regulator; see \cite{Davenport2000book}. Thus, the leading constant in our lower bound reflects the normalized class number in the imaginary quadratic case, and the normalized product of the class number and regulator in the real quadratic case. For the distribution and extreme values of $L(1,\chi_D)$, see \cite{Granville2003GeomFunctAnal}. For instance, when \(\mbf=\mbq(\mmi)\) and \(\mbf=\mbq(\sqrt{-3})\), we have
\[\kappa_{\mbq(\mmi)}=\frac{\pi}{4}\quad \text{and}\quad \kappa_{\mbq(\sqrt{-3})}=\frac{\pi}{3\sqrt{3}},\]
respectively.

For a cyclotomic field \(\mbf=\mbq(\omega_m)\), where $\omega_m$ is a primitive $m$-th root of unity and $m\ge 3$, let \(\chi^\ast\) denote the primitive Dirichlet character inducing a character \(\chi\) modulo \(m\). By the factorization of the Dedekind zeta function, we have
\[\kappa_{\mbq(\omega_m)} = \prod_{\substack{\chi\,(\bmod m)\\\chi\neq\chi_0}}L(1,\chi^\ast).\]
In particular, when $m=q$ is prime, every non-principal character modulo $q$ is primitive. Thus, Theorem \ref{thm1} gives the corresponding lower bound with an explicit leading constant for cyclotomic fields.

\section{Preliminary lemmas}
\label{sec-lem}

In this section, we present several lemmas. We begin with a classical estimate due to Landau \cite{Landau1949book}.
\begin{lem}
\label{lem-AF}
    Let \(\mbf/\mbq\) be a fixed number field of degree \(d\), and let \(A_\mbf(x)\) denote the number of nonzero integral ideals of \(\mbf\) with norm not exceeding \(x\). Then for \(x \ge 2\), we have
    \[A_\mbf(x) := \sum_{n \le x}\af(n) = \kf x+O_\mbf\bg(x^{1-\delta} \bg),\]
    where \(\delta = \delta_\mbf=2/(d+1)\).
\end{lem}
\begin{proof}
    This follows from \cite[Satz 210]{Landau1949book}.
\end{proof}

The following lemma gives an asymptotic formula, which approximates derivatives of \(\zf(s)\) on the \(1\)-line by a truncated sum.
\begin{lem}
    \label{lem-approximate}
    Let \(\mbf/\mbq\) be a fixed number field of degree \(d\). When \(T\) is sufficiently large, uniformly for all positive integers \(\ell \le (\lt T)/(\lf T)\) and \(T \le t \le 2T\), we have
    \[(-1)^\ell\zf^{(\ell)}(1+\mit) = \sum_{n \le T^{d+1}}\frac{\af(n)(\log n)^\ell}{n^{1+\mit}}+O_\mbf\bg(T^{-1+o(1)}\bg).\]
\end{lem}
\begin{proof}
    Let \(E_\mbf(x):= A_\mbf(x)-\kf x\). By Lemma \ref{lem-AF}, we have
    \begin{align}
        \label{lem-app-EF}
        E_\mbf(x) \ll_\mbf x^{1-\delta},
    \end{align}
    where \(\delta=2/(d+1)\). For \(\re(s)>1\), partial summation gives
    \begin{align}
        \label{lem-app-Abel}
        \zf(s)=\sum_{n \le X}\frac{\af(n)}{n^s}+\frac{\kf X^{1-s}}{s-1}-\frac{E_\mbf(X)}{X^s}+ s\int_X^\infty\frac{E_\mbf(u)}{u^{s+1}}\mmd u,
    \end{align}
    where \(X=T^\alpha\) with \(\alpha>0\) to be chosen. In fact, by \eqref{lem-app-EF}, the integral on the right-hand side of \eqref{lem-app-Abel} converges absolutely for \(\re(s)>1-\delta\). Hence, by analytic continuation, \eqref{lem-app-Abel} holds for \(\re(s)>1-\delta\), apart from the pole at \(s=1\). Differentiating \eqref{lem-app-Abel} \(\ell\) times and multiplying by \((-1)^\ell\) gives
    \begin{align}
        \label{lem-app-diff}
        (-1)^\ell\zf^{(\ell)}(s)&=\sum_{n \le X}\frac{\af(n)(\log n)^\ell}{n^s}+\kf X^{1-s}\sum_{j=0}^\ell\binom{\ell}{j}\frac{j!(\log X)^{\ell-j}}{(s-1)^{j+1}} - \frac{E_\mbf(X)(\log X)^\ell}{X^s} \nonumber \\
        &\quad +s\int_X^\infty \frac{E_\mbf(u)(\log u)^\ell}{u^{s+1}}\mmd u -\ell\int_X^\infty \frac{E_\mbf(u)(\log u)^{\ell-1}}{u^{s+1}}\mmd u.
    \end{align}

Now let \(s=1+\mit\), where \(T \le t \le 2T\). We estimate the last four terms on the right-hand side of \eqref{lem-app-diff}, respectively. Since
\[\binom{\ell}{j}j!=\frac{\ell!}{(\ell-j)!} \le \ell^j,\]
we have
\begin{align}
    \label{lem-app-error1}
    \Bg|\kf X^{-\mit}\sum_{j=0}^\ell\binom{\ell}{j}\frac{j!(\log X)^{\ell-j}}{(\mit)^{j+1}} \Bg|\le \frac{\kf (\log X)^\ell}{T}\sum_{j=0}^\ell \Bg(\frac{\ell}{T\log X}\Bg)^j \ll_\mbf \frac{(\log X)^\ell}{T}.
\end{align}
In the last step, we use \(\ell \le (\lt T)/(\lf T) = o(\log X)\). Moreover, Lemma \ref{lem-AF} gives
\begin{align}
    \label{lem-app-error2}
    \Bg|\frac{E_\mbf(X)(\log X)^\ell}{X^{1+\mit}} \Bg| \ll_\mbf X^{-\delta}(\log X)^\ell.
\end{align}
To estimate the integrals, observe that, for \(0 \le j \le \ell\) and \(\ell=o(\log X)\),
\begin{align*}
    \int_X^\infty \frac{(\log u)^j}{u^{1+\delta}}\mmd u = X^{-\delta}(\log X)^j \int_0^\infty\Bg(1+\frac{v}{\log X}\Bg)^j\mme^{-\delta v}\mmd v \ll_\mbf X^{-\delta}(\log X)^j.
\end{align*}
Consequently, we have
\begin{align}
    \label{lem-app-error3}
    \Bg|s\int_X^\infty \frac{E_\mbf(u)(\log u)^\ell}{u^{s+1}}\mmd u\Bg| \ll_\mbf TX^{-\delta}(\log X)^\ell,
\end{align}
and
\begin{align}
    \label{lem-app-error4}
    \Bg|\ell\int_X^\infty \frac{E_\mbf(u)(\log u)^{\ell-1}}{u^{s+1}}\mmd u\Bg|\ll_\mbf \ell X^{-\delta}(\log X)^{\ell-1}\ll_\mbf X^{-\delta}(\log X)^\ell.
\end{align}
Combining \eqref{lem-app-error1}, \eqref{lem-app-error2}, \eqref{lem-app-error3}, \eqref{lem-app-error4} with \eqref{lem-app-diff}, we get
\begin{align}
    \label{lem-app-final}
    (-1)^\ell\zf^{(\ell)}(1+\mit)=\sum_{n \le X}\frac{\af(n)(\log n)^\ell}{n^{1+\mit}}+O_\mbf\bg((\log X)^\ell\bg(T\inv +TX^{-\delta}\bg)\bg).
\end{align}

Set \(\alpha=d+1\) so that \(TX^{-\delta}=T\inv\). On the other hand, it follows from
\[\log(\log X)^\ell = \ell\log((d+1)\log T) \ll_\mbf \frac{\lo T\lt T}{\lf T}=o(\log T)\]
that \((\log X)^\ell = T^{o(1)}\) holds. Substituting these estimates into \eqref{lem-app-final} yields that
\[(-1)^\ell\zf^{(\ell)}(1+\mit) = \sum_{n \le T^{d+1}}\frac{\af(n)(\log n)^\ell}{n^{1+\mit}}+O_\mbf\bg(T^{-1+o(1)}\bg)\]
holds uniformly for all positive integers \(\ell \le (\lt T)/(\lf T)\) and \(T \le t \le 2T\). This completes the proof.    
\end{proof}

The following result provides a basic estimate, which will be used repeatedly in the subsequent proof.
\begin{lem}
    \label{lem-basic}
    Let $A,C>0$ be fixed. As $T\to\infty$, we have 
    \[\sup_{0\le j\le \frac{\lt T}{\lf T}+1}\frac{(\log_3T)^A Y_j}{\log_2T}=o(1)\quad \text{and} \quad\sup_{0\le j\le \frac{\log_3T}{\log_4T}+1}\frac{(\log_3T)^A C^j j!}{\log_2T}=o(1),\]
    where both supremum are taken over nonnegative integers $j$.
\end{lem}
\begin{proof}
    By \cite[Lemma 2.6]{Hildebrand1993JTNBordeaux}, the Laplace transform of the Dickman function is given by
    \[\widehat{\rho}(s)=\int_0^\infty \rho(u)\mme^{-su}\mmd u =\exp\Bg(\gamma+\int_0^s\frac{\mme^{-z}-1}{z}\mmd z\Bg),\quad s \in \mbc.\]
    Taking $s=-1$, we obtain \(\int_0^\infty\mme^u \rho(u)\mmd u < \infty\). Since $u^j\le j!\mme^u$ for $u\ge0$, it follows that \(Y_j \ll j!\). On the other hand, by Stirling's formula, uniformly for \(1 \le j \le (\lt T)/(\lf T)+1\), we have
    \begin{align*}
        \log\frac{(\lt T)^A C^j j!}{\lo T}
&\le A\lf T+j\log j-\lt T+O_C(j)\\
&\le -\frac{\log_3 T\log_5 T}{\log_4 T} +O_{A,C}\Bg(\frac{\log_3 T}{\log_4 T}+\lf T\Bg) \to -\infty.
    \end{align*}
    The case $j=0$ is immediate. Both assertions follow.
\end{proof}

The following lemma provides uniform estimates for the weighted Dickman function and its tail integral.
\begin{lem}
    \label{lem-tail}
    Suppose that $U\sim\log_3T$ as $T\to\infty$. Then uniformly for all positive integers \(\ell \le (\lt T)/(\lf T)\), we have \(U^\ell\rho(U)=o(1)\), and
    \[\int_U^\infty u^\ell \rho(u)\mmd u=o(1).\]
\end{lem}
\begin{proof}
    de Bruijn \cite{deBruijn1951JIndianMS} gives that
    \[\log \rho(u) = -u \Bg(\log u+\lo(u+2)-1+O\Bg(\frac{\lo(u+2)}{\log (u+2)}\Bg) \Bg)\]
    as \(u\to \infty\). Hence, for all sufficiently large \(u\), we have \(\rho(u) \le \exp((-u\log u)/2)\). Since $U\sim\lt T$ and $\ell\le\lt T/\lf T$, we have $\ell\le U/4$ for sufficiently large $T$. Hence, for $u\ge U$,
    \[u^\ell\rho(u) \le \exp\Bg(\Bg(\ell-\frac{u}{2}\Bg)\log u\Bg) \le \exp\Bg(-\frac{1}{4}u\log u\Bg).\]
    Taking \(u=U\) gives \(U^\ell\rho(U)=o(1)\). Moreover, 
    \[\int_U^\infty u^\ell \rho(u)\mmd u \le\int_U^\infty\exp\Bg(-\frac{1}{4}u\log u\Bg)\mmd u=o(1).\]
    Both estimates hold uniformly for \(\ell \le (\lt T)/(\lf T)\), completing the proof.
\end{proof}

\section{Weighted friable sums}
\label{sec-sums}

For \(x \ge 1\) and \(y \ge 2\), we define the weighted friable sum
\[
\Psi_\mbf(x,y) := \sum_{\substack{n \le x \\ P^+(n) \le y}}\af(n),
\]
where \(P^+(n)\) is the largest prime factor of \(n\), with \(P^+(1)=1\). To establish an asymptotic formula for \(\Psi_\mbf(x,y)\), we first verify that \(\af\) belongs to the class of arithmetic functions considered by de la Bret\`{e}che and Tenenbaum \cite{delaBreteche2024FACM}.

\begin{lem}
    \label{lem-aH}
    Let \(\mbf/\mbq\) be a fixed number field of degree \(d\). There exists a constant \(c_F>0\) such that
    \[\af \in \mch\Bg(1,d+1;\frac{1}{3},c_\mbf,\frac{1}{12}\Bg),\]
    where \(\mch(z,\kappa;\beta,c,\delta)\) denotes the class of arithmetic functions introduced in \cite{delaBreteche2024FACM}.
\end{lem}
\begin{proof}
    Write 
\[\zf(s)=\zeta(s)\mcb_\mbf(s),\quad \mcb_\mbf(s) = \frac{\zf(s)}{\zeta(s)}.\]
Since both \(\zf(s)\) and \(\zeta(s)\) have simple poles at \(s=1\), we have \(\mcb_\mbf(1)=\kf\). We now verify the conditions required in \cite{delaBreteche2024FACM}. First, we take
\[\beta=\frac{1}{3},\quad \delta = \frac{1}{12},\quad \beta+\delta = \frac{5}{12}<\frac{3}{5}.\]
By the classical zero-free region for \(\zeta(s)\) (see \cite[Theorem 3.8]{Titchmarsh1986book}), the function \(\mcb_\mbf(s)\) admits a holomorphic continuation to 
\[\sigma>1-\frac{c_\mbf}{(1+\log^+|t|)^{7/5}}\]
for some sufficiently small constant \(c_\mbf>0\). Moreover, by \cite[Theorem 3.11]{Titchmarsh1986book}, we have
\begin{align}
    \label{lem-aH-Tom}
    \frac{1}{\zeta(s)}\ll \log|t|
\end{align}
in the corresponding zero-free region for sufficiently large \(|t|\). On the other hand, Scourfield \cite[Lemma 2.6]{Scourfield2004JTNBordeaux} gives
\begin{align}
    \label{lem-aH-Jerry}
    \zf(s) \ll_\mbf (\log|t|)^{2/3}\lo |t|
\end{align}
in a region containing the one considered above. Consequently, by \eqref{lem-aH-Tom} and \eqref{lem-aH-Jerry},
\[\mcb_\mbf(s)=\frac{\zf(s)}{\zeta(s)}\ll_\mbf(\log|t|)^{5/3}\lo |t|\ll_\mbf(1+|t|)^{11/12},\]
uniformly in this region. Thus, the required analytic condition is satisfied.

Next, write                                                                \[\mcb_\mbf(s)=\sum_{n=1}^\infty \frac{b_\mbf(n)}{n^s},\]                      
where \(b_\mbf=\af*\mu\). For \(y \ge 2\), set
\[\mcb_\mbf(s;y):=\sum_{\substack{n \ge 1 \\ P^+(n) \le y}}\frac{b_\mbf(n)}{n^s}\]
By the Euler product,
\[
\mcb_\mbf(s;y)=\prod_{p\le y}\bg(1-p^{-s}\bg)\prod_{p\le y}\prod_{\mathfrak p\mid p}\bg(1-(N\mathfrak p)^{-s}\bg)^{-1}.
\]
To estimate this product, for \(\mbe=\mbf,\mbq\), define
\[
Z_\mbe(s;y)=\prod_{N\mathfrak p\le y}\bg(1-(N\mathfrak p)^{-s}\bg)^{-1},
\]
Fix a sufficiently small \(\eta \in (0,1/60)\). By \cite[Lemma 4.1]{Scourfield2004JTNBordeaux}, for \(\mbe=\mbf,\mbq\), we have
\[Z_\mbe(s;y)=\zeta_\mbe(s)(s-1)\log y\widehat{\rho}((s-1)\log y)\bg(1+O_\mbe\bg(\mme^{-(\log y)^{3/5-\eta}}\bg)\bg)\]
uniformly for
\[
\sigma\ge1-(\log y)^{-2/5-\eta},\quad |t|\le \mme^{(\log y)^{3/5-\eta}},
\]
where \(\widehat{\rho}\) denotes the Laplace transform of the Dickman function. Then, taking the quotient gives
\begin{align}
    \label{lem-aH-quotient}
    \frac{Z_\mbf(s;y)}{Z_\mbq(s;y)} = \mcb_\mbf(s)\bg(1+O_\mbf\bg(\mme^{-(\log y)^{3/5-\eta}}\bg)\bg).
\end{align}
On the other hand, 
\[\mcb_\mbf(s;y)= \frac{Z_\mbf(s;y)}{Z_\mbq(s;y)}\prod_{\substack{p \le y,\mathfrak{p}\mid p\\N\mathfrak{p}>y}}\bg(1-\bg(N\mathfrak p\bg)^{-s}\bg)\inv.\]
For \(\sigma \ge 3/4\), we have
\[\log\prod_{\substack{p \le y,\mathfrak{p}\mid p\\N\mathfrak{p}>y}}\bg(1-\bg(N\mathfrak p\bg)^{-s}\bg)\inv \ll_\mbf \sum_{f=2}^d\sum_{p>y^{1/f}}p^{-3f/4}\ll_\mbf y^{-1/4}.\]
Consequently, combining this with \eqref{lem-aH-quotient}, we obtain
\[\mcb_\mbf(s;y)=\mcb_\mbf(s)\bg(1+O_\mbf\bg(\mme^{-(\log y)^{3/5-\eta}}+y^{-1/4}\bg)\bg).\]
Since \(5/12+\eta<3/5\), we have
\[\mcb_\mbf(s;y)=\mcb_\mbf(s)+O_\mbf\bg(\mme^{-(\log y)^{5/12}}\bg)\]
uniformly for the above region. This establishes the required approximation for the truncated Dirichlet series.

Finally, since \(0\le \af(n)\le \tau_d(n)\le\tau_{d+1}(n)\) and
\[\sum_{n=1}^\infty\frac{\tau_{d+1}(n)}{n^s}=\zeta(s)^{d+1},\]
we may take
\[
\af^\dagger(n)=\tau_{d+1}(n)
\]
as a majorant, where \(\tau_k(n)\) denotes the \(k\)-fold divisor function. The majorant condition is therefore satisfied. This completes the proof.  
\end{proof}

The following two propositions provide the estimates for weighted friable sums needed in the proof of Theorem \ref{thm1}.

\begin{pro}
\label{pro1}
    Let \(\mbf/\mbq\) be a fixed number field, and let \(u=\log x/\log y\). Then, uniformly for \(x \ge y \ge \exp\bg((\log x)^{2/3}\bg)\), we have
    \[\Psi_\mbf(x,y)=\kf x \rho(u)+O_\mbf\Bg(x\rho(u)\frac{\log(2u)}{\log y}\Bg).\]
\end{pro}
\begin{proof}
    By Lemma \ref{lem-aH}, we apply \cite[Corollary 1.3 \((\mathrm{\romannumeral 2})\)]{delaBreteche2024FACM} with \(z=1\). Following the notation of \cite{delaBreteche2024FACM}, we have \(\psi_1(u)=\rho(u)\) and \(R_1(u)\asymp\rho(u)\). Moreover,    \[a_0(\af)=\lim_{w \to 0}\frac{w\zf(1+w)}{1+w}=\kf.\]
    Since \(\rho(u)\asymp1\) for \(1\le u\le2\), the additional error term in \cite[Corollary 1.3 \((\mathrm{\romannumeral 2})\)]{delaBreteche2024FACM}, supported near \(u=1\), can be absorbed into the stated error term. This completes the proof. 
\end{proof}

Let
\[b:=\lfloor (\lo T)^3 \rfloor,\quad y:= \frac{\log T}{40b},\quad L:=\exp(\lo T \lt T).\]
With these choices, we obtain the following asymptotic formula.

\begin{pro}
    \label{pro2}
    Let \(\mbf/\mbq\) be a fixed number field. As \(T \to \infty\), uniformly for all positive integers $\ell\le (\lt T)/(\lf T)$, we have
    \[\sum_{\substack{n \le L \\ P^+(n) \le y}}\frac{\af(n)(\log n)^\ell}{n} = \bg(\kf Y_\ell+o_\mbf(1)\bg)(\lo T)^{\ell+1}.\]
\end{pro}
\begin{proof}
    By partial summation,
    \begin{align}
        \label{pro2-Abel}
        S &:=S_{\mbf,\ell}(L,y)= \sum_{\substack{n \le L \\ P^+(n) \le y}}\frac{\af(n)(\log n)^\ell}{n} \nonumber \\
        & = \frac{(\log L)^\ell}{L}\Psi_\mbf(L,y) + \int_1^L\Psi_\mbf(x,y)\frac{(\log x)^\ell-\ell(\log x)^{\ell-1}}{x^2}\mmd x. 
    \end{align}
    When \(1 \le x \le y\), by Lemma \ref{lem-AF}, we have
    \[\Psi_\mbf(x,y)=A_\mbf(x)=\kf x +O_\mbf\bg(x^{1-\frac{2}{d+1}}\bg).\]
    On the other hand, when \(y \le x \le L\), since \((\log L)^{2/3}=o(\lo  T)\), Proposition \ref{pro1} applies uniformly for this range. Thus, for \(1 \le x \le L\), 
    \begin{align}
        \label{pro2-Psif}
        \Psi_\mbf(x,y)=\kf x\rho\Bg(\frac{\log x}{\log y}\Bg)+\mce_\mbf(x),
    \end{align}
    where
    \begin{equation*}
    \mce_\mbf(x) \ll_\mbf
    \begin{cases}
        x^{1-\frac{2}{d+1}}, ~ &\operatorname{if} 1 \le x \le y, \\
		x\rho\Bg(\frac{\log x}{\log y}\Bg)\frac{\log(2\log x/\log y)}{\log y}, ~ &\operatorname{if} y \le x \le L.
    \end{cases}
\end{equation*}

Substituting \eqref{pro2-Psif} into \eqref{pro2-Abel} yields that
\begin{align}
    \label{pro2-S}
    S=\kf(\log L)^\ell\rho\Bg(\frac{\log L}{\log y}\Bg)+\kf \int_1^L\rho\Bg(\frac{\log x}{\log y}\Bg)\frac{(\log x)^\ell-\ell(\log x)^{\ell-1}}{x}\mmd x + \mce,
\end{align}
where
\begin{align*}
    \mce &= O_\mbf\Bg(C_\mbf^\ell \ell!+ \frac{(\log L)^\ell}{\log y} \rho\Bg(\frac{\log L}{\log y}\Bg)\log\Bg(\frac{2\log L}{\log y}\Bg)\Bg)\\
    &\quad +O_\mbf\Bg(\frac{1}{\log y}\int_y^L\rho\Bg(\frac{\log x}{\log y}\Bg)\log\Bg(\frac{2\log x}{\log y}\Bg)\frac{(\log x)^\ell+\ell(\log x)^{\ell-1}}{x}\mmd x \Bg),
\end{align*}
and \(C_\mbf>0\) depends only on \(\mbf\). Making the change of variables \(u=(\log x)/(\log y)\), the first two terms on the right-hand side of the above identity equal
\[\kf(\log y)^\ell U^\ell\rho(U) + \kappa_F(\log y)^{\ell+1}
\int_0^U u^\ell\rho(u)\mmd u-\kappa_F\ell(\log y)^\ell\int_0^U u^{\ell-1}\rho(u)\mmd u,\]
where \(U=(\log L)/(\log y)\). Note that \(U=(1+o(1))\lt T\), by using Lemma \ref{lem-basic} and Lemma \ref{lem-tail}, this equals
\begin{align}
    \label{pro2-main}
    \bg(\kf Y_\ell +o_\mbf(1)\bg)(\log y)^{\ell+1}.
\end{align}

It remains to estimate the error term \(\mce\). In fact, it satisfies
\[\frac{|\mce|}{(\log y)^{\ell+1}}\ll_\mbf \frac{C_\mbf^\ell \ell!}{(\log y)^{\ell+1}}+\frac{Y_{\ell}+Y_{\ell+1}}{\log y}+\frac{\ell(Y_{\ell-1}+Y_\ell)}{(\log y)^2}+\frac{U^\ell\rho(U)\log(2U)}{(\log y)^2}.\]
By repeatedly applying Lemma \ref{lem-basic}, one deduces, for large \(y\),
\begin{align}
    \label{pro2-error}
    \mce = o\bg((\log y)^{\ell+1}\bg).
\end{align}
By the definitions of \(b\) and \(y\), we have \(\log y=\log_2 T-3\log_3 T-\log 40+o(1)\). Consequently, uniformly for \(\ell\le (\log_3 T)/(\log_4 T)\),
\[\Bg|\Bg(\frac{\log y}{\log_2 T}\Bg)^{\ell+1}-1\Bg|\ll\frac{(\ell+1)\log_3 T}{\log_2 T}.\]
It follows from Lemma \ref{lem-basic}, with $A=2$, that
\[Y_\ell\Bg|\Bg(\frac{\log y}{\log_2 T}\Bg)^{\ell+1}-1\Bg|\ll
\frac{Y_\ell(\ell+1)\log_3 T}{\log_2 T}=o(1).\]
Therefore, combining \eqref{pro2-main} and \eqref{pro2-error}, we obtain
\[S=\bg(\kf Y_\ell +o_\mbf(1)\bg)(\lo T)^{\ell+1}\]
uniformly for all positive integers $\ell\le (\lt T)/(\lf T)$. This completes the proof.
\end{proof}

\section{Proof of Theorem \ref{thm1}}
\label{sec-pthm1}
In this section, we prove Theorem \ref{thm1} using an approach similar to Yang \cite[Section 4]{Ddyang2022Mathematika}. Let $N=T^{1/4}$ and $X=T^{d+1}$. For a function $r:\mbn\to\{0,1\}$, define the resonator
\[\mcr(t)=\sum_{n\le N}r(n)n^{-\mit}.\] 
As in \cite{Soundararajan2008MathAnn}, let \(\Phi\) be a smooth function compactly supported in \([1,2]\) to satisfy \(0 \le \Phi(t) \le 1\) always and \(\Phi(t)=1\) for \(5/4 < t < 7/4\). Repeated integration by parts shows that its Fourier transform \(\widetilde{\Phi}(\xi) \ll_v |\xi|^{-v}\) holds for every positive integers \(v\).

Consider the integrals
\[I_0:=I_0(\mcr,T) = \int_\mbr |\mcr(t)|^2\Phi\Bg(\frac{t}{T}\Bg)\mdt,\]
and
\[I_1:=I_1(\mcr,T)=\int_\mbr (-1)^\ell\zf^{(\ell)}(1+\mit)|\mcr(t)|^2\Phi\Bg(\frac{t}{T}\Bg)\mdt.\]
Clearly, we obtain
\begin{align}
    \label{max1}
    \max_{t \in [T,2T]}\bg|\zf^{(\ell)}(1+\mit)\bg| \ge \frac{|I_1|}{I_0}.
\end{align}
Following the argument in \cite[p. 471]{Soundararajan2008MathAnn}, the Cauchy–Schwarz inequality and the rapid decay of \(\widetilde{\Phi}\) give
\begin{align}
    \label{I0-equality}
    I_0 = T \widetilde{\Phi}(0)\bg(1+O\bg(T\inv\bg)\bg)\sum_{n \le N}|r(n)|^2.
\end{align}

On the other hand, Lemma \ref{lem-approximate} gives
\[I_1 = \int_\mbr \sum_{k \le X}\frac{\af(k)(\log k)^\ell}{k^{1+\mit}}|\mcr(t)|^2\Phi\Bg(\frac{t}{T}\Bg)\mdt+O_\mbf\bg(T^{-1+o(1)}I_0\bg).\]
Expanding \(|\mcr(t)|^2\), then we obtain
\[I_1 = T\sum_{k \le X}\frac{\af(k)(\log k)^\ell}{k}\sum_{m,n\le N}r(m)\overline{r(n)}\widetilde{\Phi}\Bg(T\log \frac{km}{n}\Bg)+O_\mbf\bg(T^{-1+o(1)}I_0\bg).\]
For $km\neq n$, since $n\le N$ and $km,n$ are distinct positive integers,
\[\Bg|\log\frac{km}{n}\Bg| \ge \frac{1}{N+1}.\]
Moreover, using $\af(k)\le\tau_d(k)$, we obtain
\[\sum_{k \le X}\frac{\af(k)(\log k)^\ell}{k} \le (\log X)^\ell \sum_{k \le X}\frac{\tau_d(k)}{k} \ll_\mbf (\log X)^{\ell+d}.\]
Here, we use the classical estimate for the divisor function, namely
\[\sum_{k \le X}\frac{\tau_d(k)}{k} \ll_d (\log X)^d.\]
Since \(X=T^{d+1}\) and \(\ell \le (\lt T)/(\lf T)\), it follows that
\[\sum_{k \le X}\frac{\af(k)(\log k)^\ell}{k} = T^{o(1)}.\]
Consequently, by the rapid decay of \(\widetilde{\Phi}\), the contribution of
all terms with $km\ne n$ is negligible. More precisely, by applying the Cauchy–Schwarz inequality, we obtain
\begin{align}
    \label{I1-equality}
    I_1=T\widetilde{\Phi}(0) \sum_{km=n\le N}\frac{\af(k)(\log k)^\ell}{k}r(m)\overline{r(n)}+O_\mbf\Bg(T^{o(1)}\sum_{n\le N}|r(n)|^2 \Bg).
\end{align}
Since \(\af(k)\) and \(r(n)\) are both non-negative, by \eqref{max1}, we have
\begin{align}
    \label{max2}
    \max_{t \in [T,2T]}\bg|\zf^{(\ell)}(1+\mit)\bg| \ge \bg(1+O\bg(T\inv\bg)\bg)\mcq_\ell(r)+O_\mbf\bg(T^{-1+o(1)}\bg),
\end{align}
where
\[\mcq_\ell(r):=\mcq_\ell(r,\mbf,N) = \Bg|\sum_{km=n\le N}\frac{\af(k)(\log k)^\ell}{k}r(m)\overline{r(n)}\Bg| \Bg/\Bg(\sum_{n\le N}|r(n)|^2 \Bg).\]

We next focus on obtaining an effective lower bound for \(\mcq_\ell(r)\). To this end, using the notation introduced before Proposition \ref{pro2} in Section \ref{sec-sums}, we define
\[\mcp(y,b):=\prod_{p \le y}p^{b-1}\quad \text{and}\quad \mcm:=\{n:n \mid \mcp(y,b)\}.\]
By the prime number theorem, we have \(\mcp(y,b) \le N\). Take $r(n)$ to be the characteristic function of the set \(\mcm\). In this case,
\[\sum_{n\le N}|r(n)|^2=|\mcm|,\]
and
\begin{align}
    \label{Q-equality-1}
    \mcq_\ell(r)=\frac{1}{|\mcm|}\sum_{\substack{n \in \mcm \\ k \mid n}}\frac{\af(k)(\log k)^\ell}{k} = \sum_{k \in \mcm}\frac{\af(k)(\log k)^\ell}{k}\sum_{\substack{n \in \mcm \\ k \mid n}}\frac{1}{|\mcm|} .
\end{align}
Define \(\mck\) as the set of positive integers \(k\) satisfying \(k \le L\) and \(P^+(k) \le y\). Here \(L=\exp(\lo T \lt T)\), which is as defined before. Set
\[\mck_1=\Bg\{k \in \mck:\Omega(k) \le \frac{b}{(\lo T)^2}\Bg\},\]
where $\Omega(k)$ denotes the number of prime factors of $k$, counted
with multiplicity. It is clear that \(\mck_1\subseteq \mcm\). Furthermore, set \(\mck_2 = \mck\setminus\mck_1\). 

For any \(k \in \mck_1\), write
\[k = \prod_{p \le y}p^{\alpha_p}.\]
Since
\[\sum_{p \le y}\alpha_p=\Omega(k) \le \frac{b}{(\lo T)^2},\]
for the inner sum in the right-hand side of \eqref{Q-equality-1}, we have
\[\sum_{\substack{n \in \mcm \\ k \mid n}}\frac{1}{|\mcm|}=\prod_{p \le y}\Bg(1-\frac{\alpha_p}{b}\Bg) \ge 1-\frac{\Omega(k)}{b}\ge 1-\frac{1}{(\lo T)^2}.\]
Substituting this lower bound into \eqref{Q-equality-1} gives
\begin{align}
\label{Q-equality-2}
    \mcq_\ell(r) \ge \Bg( 1-\frac{1}{(\lo T)^2}\Bg) \sum_{k \in \mck_1}\frac{\af(k)(\log k)^\ell}{k}.
\end{align}
We now treat the contribution from the set \(\mck_2\). To this end, fix \(\varepsilon\in(0,\log2)\). By Rankin's trick and $\af(k)\le\tau_d(k)$, we have
\begin{align*}
    \sum_{k \in \mck_2} \frac{\af(k)(\log k)^\ell}{k} &\le (\log L)^\ell \mme^{-\varepsilon b/(\lo T)^2}\sum_{P^+(k)\le y} \frac{\af(k)\mme^{\varepsilon\Omega(k)}}{k} \\
    & \le (\log L)^\ell \mme^{-\varepsilon b/(\lo T)^2} \prod_{p \le y}\Bg(1-\frac{\mme^\varepsilon}{p}\Bg)^{-d}.
\end{align*}
Then, by using Mertens' theorem, 
\[\sum_{k \in \mck_2} \frac{\af(k)(\log k)^\ell}{k} \ll_{\mbf,\varepsilon} (\log L)^\ell \mme^{-\varepsilon b/(\lo T)^2}(\log y)^{d\mme^\varepsilon}.\]
Since \(\lo y = O(\lt T)=o(\lo T)\),
\[\frac{b}{(\lo T)^2}=\bg(1+o(1)\bg)\lo T \quad \text{and} \quad \ell\lo L \ll \frac{(\lt T)^2}{\lf T}=o(\lo T),\]
we obtain
\begin{align}
    \label{K2-bound}
    \sum_{k \in \mck_2} \frac{\af(k)(\log k)^\ell}{k} \ll_{\mbf,\varepsilon} \exp(-\varepsilon \lo T+o(\lo T) ) = o_\mbf(1).
\end{align}
By the definition of \(\mck\) and the construction of \(\mck_2\), \eqref{K2-bound} yields 
\[\sum_{k \in \mck_1} \frac{\af(k)(\log k)^\ell}{k} =\sum_{\substack{k\le L\\P^+(k)\le y}}\frac{\af(k)(\log k)^\ell}{k}+o_\mbf(1). \]
Furthermore, by using Proposition \ref{pro2}, we obtain
\begin{align}
    \label{K1-bound}
    \sum_{k \in \mck_1} \frac{\af(k)(\log k)^\ell}{k} =\sum_{\substack{k\le L\\P^+(k)\le y}}\frac{\af(k)(\log k)^\ell}{k}+o_\mbf(1) = \bg(\kf Y_\ell + o_\mbf(1)\bg) (\lo T)^{\ell+1}.
\end{align}
Moreover, Lemma \ref{lem-basic} gives 
\[\frac{Y_\ell}{(\lo T)^2} = o(1)\]
uniformly for $\ell\le(\log_3T)/(\log_4T)$. Thus, combining this estimate with \eqref{Q-equality-2} and \eqref{K1-bound} gives
\[\mcq_\ell(r) \ge \bg(\kf Y_\ell + o_\mbf(1)\bg) (\lo T)^{\ell+1}.\]
Finally, substituting this into \eqref{max2} completes the proof of Theorem \ref{thm1}.

\section*{Acknowledgments}
Zikang Dong was supported by the National Natural Science Foundation of China (Grant No. 1240011770). Qiyu Yang was supported by the Natural Science Foundation of Henan Province (Grant No. 252300421782) and the National Natural Science Foundation of China (Grant No. 12601011). Guang-Liang Zhou was supported by the National Natural Science Foundation of China (Grant No. 12401009).
	
	\bibliographystyle{siam}
    \bibliography{references}
\end{document}